\documentclass[12pt,twoside]{article}
\usepackage{amssymb,amsmath,latexsym}

\begin{document}

\title{Decay bounds on eigenfunctions and the singular spectrum of unbounded Jacobi matrices}

\author{Jan Janas, Serguei Naboko and G\"unter Stolz}

\date{ }

\newtheorem{theo}{Theorem}[section]
\newtheorem{defini}[theo]{Definition}
\newtheorem{proposi}[theo]{Proposition}
\newtheorem{lemma}[theo]{Lemma}
\newtheorem{coro}[theo]{Corollary}
\newtheorem{rem}[theo]{Remark}
\newtheorem{ex}[theo]{Example}
\newtheorem{conj}[theo]{Conjecture}
\newtheorem{proof}[theo]{Proof}
\newcommand{\CC}{{\mathbb C}}
\newcommand{\NN}{{\mathbb N}}
\newcommand{\RR}{{\mathbb R}}
\newcommand{\SM}{{\mathbb S}}
\newcommand{\TM}{{\mathbb T}}
\newcommand{\ZZ}{{\mathbb Z}}

\newcommand{\HH}{{\bf H}}
\newcommand{\Pp}{{\cal P}}
\newcommand{\Zz}{{\cal Z}}
\newcommand{\PP}{{\bf P}}
\newcommand{\EE}{{\bf E}}
\newcommand{\Bb}{{\cal B}}
\newcommand{\Ee}{{\cal E}}
\newcommand{\Ff}{{\cal F}}
\newcommand{\Ww}{{\cal W}}
\newcommand{\Ss}{{\cal S}}
\newcommand{\Oo}{{\cal O}}
\newcommand{\Tr}{\mbox{\rm Tr}}
\newcommand{\Rr}{{\cal R}}
\newcommand{\Nn}{{\cal N}}
\newcommand{\Cc}{{\cal C}}
\newcommand{\Jj}{{\cal J}}
\newcommand{\Ll}{{\cal L}}
\newcommand{\dyn}{{\cal S}_{\epsilon,\omega}}
\newcommand{\dynd}{{\cal S}_{\delta,\omega}}
\newcommand{\dynl}{{\cal S}_{\epsilon,\omega_l}}
\newcommand{\dynopm}{{\cal S}_{0,\pm}}
\newcommand{\dyno}{{\cal S}_{0,\omega}}
\newcommand{\dynpm}{{\cal S}_{\epsilon,\pm}}
\newcommand{\cor}{{\mbox{\rm\tiny cor}}}
\newcommand{\eig}{{\mbox{\rm\tiny eig}}}
\newcommand{\dis}{{\mbox{\rm\tiny dis}}}
\newcommand{\dev}{{\mbox{\rm\tiny dev}}}
\def\essinf{\mathop{\rm ess\,inf}}
\def\esssup{\mathop{\rm ess\,sup}}

\maketitle

\begin{abstract}

Bounds on the exponential decay of generalized eigenfunctions of
bounded and unbounded selfadjoint Jacobi matrices in $\ell^2(\NN)$
are established. Two cases are considered separately and lead to
different results: (i) the case in which the spectral parameter lies
in a general gap of the spectrum of the Jacobi matrix and (ii) the
case of a lower semi-bounded Jacobi matrix with values of the
spectral parameter below the spectrum. It is demonstrated by
examples that both results are sharp.

We apply these results to obtain a ``many barriers-type'' criterion for the existence of square-summable generalized eigenfunctions of an unbounded Jacobi matrix at almost every value of the spectral parameter in suitable open sets. In particular, this leads to examples of unbounded Jacobi matrices with a spectral mobility edge, i.e.\ a transition from purely absolutely continuous spectrum to dense pure point spectrum.

\end{abstract}

\setcounter{equation}{0}
\section{Introduction} \label{sec1}

One of the central tools in the spectral theory of differential and finite difference operators, in particular Schr\"odinger operators and their discrete counterparts, are results on the asymptotic behavior of generalized eigenfunctions. Here we are concerned with establishing such results for unbounded Jacobi matrices and relating them to spectral properties of the associated self-adjoint operators.

In the first part of the paper we consider general unbounded self-adjoint Jacobi
matrices $J$ on $\ell^2(\NN)$ and assume that $\lambda \in \RR$ lies
in a spectral gap of $J$. We will use a ``discrete'' and rather
simple version of a technique introduced in \cite{BCH} to prove
upper bounds on the exponential decay of generalized eigenfunctions
of $J$ to $\lambda$. The decay bound for eigenfunctions of
Schr\"odinger operators found in \cite{BCH} improved on longstanding
bounds obtained through the so-called Combes-Thomas method
\cite{CT}. Our results on exponential decay are stated in Section~\ref{sec2} and proven in Section~\ref{sec3}.

While these results are quite general and their proof, due to the
discrete one-dimensional setting, quite elementary, the obtained
bounds are remarkably sharp in several respects. This will be
understood in Section~\ref{sec4}, where we will consider
two concrete classes of unbounded Jacobi matrices for which the
exact asymptotics of generalized eigenfunctions can be obtained.

Combes-Thomas type estimates like the ones proven here are a
frequently used tool in the spectral analysis of differential and
finite difference operators. Some new applications are given in Sections~\ref{sec5} and \ref{sec6} below. We also mention a recent paper by J.\ Breuer \cite{Breuer}, where a
similar but weaker estimate for the matrix elements of $(J-z)^{-1}$
was found and applied to study spectral and dynamical properties of
random Jacobi matrices.

Our main application of the results on exponential decay is a criterion for the existence of $\ell^2$-solutions of $Ju=\lambda u$ at (Lebesgue) almost every energy $\lambda$ in an interval $I$. This describes types of Jacobi matrices $J$ which coincide on infinitely many growing intervals with a Jacobi matrix $J_0$ which has $I$ as a spectral gap. Physically, these intervals can be thought of as a series of barriers preventing quantum mechanical transport under the time evolution for the hamiltonian $J$. Consequences are the absence of absolutely continuous spectrum for $J$ in $I$ and, after adding generic rank-one perturbations, pure point spectrum in $I$. For discrete and continuous one-dimensional Schr\"odinger operators results of this kind were first obtained by Kirsch, Molchanov and Pastur in \cite{KMP1, KMP2}. Subsequently, their ideas have been extended to multi-dimensional Schr\"odinger operators \cite{GJMS, Stolz}, as well as to finite difference operators on strips \cite{KMPV}.

We state and prove a general barriers-type criterion for unbounded Jacobi matrices in Section~\ref{sec5}. In the final Section~\ref{sec6} we study a specific class of unbounded Jacobi matrices, which give rise to a mobility edge. Their spectrum covers the entire real line, it is purely absolutely continuous outside an interval $[-c,c]$ for some $c>0$, and dense pure point in $(-c,c)$. The latter will follow from our criterion in Section~\ref{sec5} together with a Weyl-sequence argument, while the claim on absolute continuity will be a consequence of the general results on the asymptotics of solutions of difference equations in \cite{JM}.

\section{Exponential decay bounds for the resolvent} \label{sec2}

Let $\langle \cdot, \cdot \rangle$ be the inner product in
$\ell^2(\NN)$ and denote by $e_n$ the canonical basis in
$\ell^2(\NN)$. Before we proceed further recall a general result from \cite{Shubin}, going back to \cite{Demko} and \cite{Kershaw}, which says that for any bounded operator
$A$ in $\ell^2(\NN)$ with the band matrix $a_{ij}=\langle
Ae_{j},e_{i}\rangle $ having the bounded inverse $A^{-1}$, the
entries $b_{ij}$ of $A^{-1}$ satisfy the estimates:

$$ |b_{ij}|\leq{Cr^{|i-j|}}, \:i,j\in \NN $$
for some $C>0$ and $ r\in(0,1)$. In the proof of this estimate the
boundedness of $A$ is used in an essential way.

Consider a Jacobi matrix $J$ defined by the difference relations
\begin{equation} \label{eq1.1}
(Ju)(n) = \lambda_{n-1} u(n-1) + q_n u(n) + \lambda_n u(n+1),
\quad n \in \NN,
\end{equation}
and boundary condition $u(0)=0$, or, in equivalent matrix
notation,
\begin{equation} \label{eq1.2}
J = \left( \begin{array}{cccc} q_1 & \lambda_1 & & \\ \lambda_1 &
q_2 & \lambda_2 & \\ & \lambda_2 & q_3 & \ddots \\ & & \ddots &
\ddots \end{array} \right) \;.
\end{equation}
Here, the ``discrete potential'' $q_n$ is real and the ``weights''
$\lambda_n$ are positive for $n\in \NN$. For later use we also adopt
the convention $\lambda_0 =1$. Sometimes it is also convenient to
represent $J$ in the form $J = S\Lambda + \Lambda S^* + Q $, where
$\Lambda$ (respectively $Q$) are the diagonal operators given by
$\lambda_k$ (respectively $q_k$) in the canonical basis $e_n$, $n\in
\NN$ , and $S$ is the unilateral shift $Se_n = e_{n+1}$.

Under the assumption
\begin{equation} \label{eq1.3}
\sum_{n=1}^{\infty} \frac{1}{\lambda_n} = \infty
\end{equation}
this defines a unique self-adjoint operator on $\ell^2(\NN)$,
e.g.\ \cite{Berezanskii}, which will also be denoted by $J$. \\

 The main result of the first part of this paper is an
upper bound for the matrix elements of the resolvent of $J$:


\begin{theo} \label{thm1}
Let $J^*=J$ and assume that $\liminf_{n\to\infty}{\lambda_n}
> 0 $.

(a) Let $(r,s)$ be a finite gap in the spectrum of $J$. Then there
exists a constant $\eta>0$ such that for all $\lambda \in
(r,s)$ and all $n\in \NN$,
\begin{equation} \label{eq1.4}
| \langle (J-\lambda)^{-1} e_1, e_n \rangle |
\le4\max\left\{(\lambda-r)^{-1},(s-\lambda)^{-1}\right\}\exp
\left( -\eta \sqrt{(\lambda-r)(s-\lambda)} \sum_{k=1}^{n-1}
\frac{1}{\lambda_k} \right).
\end{equation}

(b) Suppose that $J$ is bounded from below and denote $d=\inf \sigma(J)$. There exists a constant $\eta>0$ such that for all $\lambda \in (-\infty,d)$ and all $n\in \NN$,
\begin{equation} \label{eq1.4b}
| \langle (J-\lambda)^{-1} e_1, e_n \rangle | \le 4(d-\lambda)^{-1} \exp \left( -\eta \sqrt{d-\lambda} \sum_{k=1}^{n-1} \frac{1}{\lambda_k} \right).
\end{equation}
The corresponding result holds if $J$ is bounded above and $\lambda \in (\sup \sigma(J),\infty)$.

\end{theo}

In the case $\lim_{n\to\infty}{\lambda_n} = +\infty $ the above
estimates can be slightly improved. This is the content of the next two results.

\begin{theo} \label{thm2}
Let $J^* = J$. Suppose that $\lim_{n\to\infty}{\lambda_n} =
+\infty $ and let $(r,s)$ be a gap in the spectrum of $J$. Then for
arbitrary $\epsilon\in(0,1/2)$ there exists $N=N(\epsilon)$ such
that
\begin{equation} \label{eq1.5}
| \langle (J-\lambda)^{-1} e_1, e_n \rangle | \le
\frac{s-r}{\epsilon(\lambda-r)(s-\lambda)}\exp \left(
-(\frac{1}{2} -\epsilon) \sqrt{(\lambda-r)(s-\lambda)}
\sum_{k=N}^{n-1} \frac{1}{\lambda_k} \right),
\end{equation}
for all $\lambda \in (r,s)$ and for all $n>N.$
\end{theo}

For a Jacobi operator $J$ which is bounded from below and
$\lambda$ below the bottom of the spectrum the estimates given in
the above results can be further improved. Indeed, we have
\begin{theo}\label{thm3}
Suppose that $J=J^* $ associated to the weights $\lambda_n$ and
the potential $q_n$ is bounded from below by $d$ and
$\lim_{n\to\infty}{\lambda_n} = +\infty $. Fix
$\epsilon\in(0,1) $ and complex $\lambda$ such that $\Re{\lambda}<d$. Then
there exists $N=N(\epsilon,\lambda)$ such that
 \begin{eqnarray}\label{eq1.6}
| \langle(J-\lambda)^{-1} e_1, e_n \rangle |
\le[(d-\Re{\lambda})\epsilon]^{-1}\exp \left( -(1-\epsilon)
\sqrt{d-\Re{\lambda}} \sum_{k=N}^{n-1} \frac{1}{\sqrt\lambda_k}
\right),
\end{eqnarray}
for $n>N.$
\end{theo}

In Section~\ref{sec4} we will demonstrate by examples that the constants in the exponents on the right hand sides of (\ref{eq1.5}) and (\ref{eq1.6}) are optimal in the sense that $1/2-\epsilon$ and $1-\epsilon$, respectively, can not be replaced by $1/2+\epsilon$ and $1+\epsilon$ for any $\epsilon>0$.

The sequence $v$ defined by $v(n) = \langle (J-\lambda)^{-1} e_1,
e_n \rangle$ is a generalized eigenfunction of $J$, meaning that
it solves (\ref{eq1.1}) for $n\ge 2$ without satisfying the
boundary condition at $0$. As $J$ is in the limit point case at
$+\infty$, it is (up to constant multiples) the unique generalized
eigenfunction which is square-summable. Thus we may understand
(\ref{eq1.1}) as a bound on the decay of generalized
eigenfunctions. As $\lambda \not\in \sigma(J)$ (the spectrum of
$J$), $v$ is not an eigenfunction, but the above results also
provide bounds on eigenfunctions for perturbations of $J$: If
$\tilde{J} = J+A$, where $A$ is a finite Jacobi matrix, and if
$u=(u(n))$ is an eigenfunction of $\tilde{J}$ for an eigenvalue
$\lambda \in (r,s)$, then $u$ satisfies the bound given by the
right hand side of (\ref{eq1.4}). This is obvious as $u(n)$ for
large $n$ coincides with the unique $\ell_2$-generalized
eigenfunction of $J$.

We point out two specific features of the exponent on the right of
(\ref{eq1.4}): It describes the large $n$ asymptotics as well as the
asymptotics as $\lambda$ approaches the spectrum of $J$. The large
$n$ asymptotics, characterized through the sum $\sum 1/\lambda_k$,
generalizes the ``Schr\"odinger case'' $\lambda_n \equiv 1$, where
generalized eigenfunctions for $\lambda$ in a spectral gap decay
exponentially in $n$. As $\lambda$ approaches a spectral edge of
$J$, that is, either $r$ or $s$, the constant in front of $\sum
1/\lambda_k$ is proportional to the {\it square root of the
distance} of $\lambda$ to the spectrum. This improvement over the
original Combes-Thomas method, which merely provides a term which is
linear in the distance, is due to the new ideas introduced in
\cite{BCH}. We will comment on this at the end of the proof of
Theorem~\ref{thm1} in Section~\ref{sec3}.

The proofs of Theorems~\ref{thm1} to \ref{thm3} allow for several generalizations. As an example, we state the following generalization of Theorem~\ref{thm1}(a), which will be used in our applications in Section~\ref{sec5}. Similar generalizations could be formulated for our other results.

\begin{theo} \label{thm4}
Under the conditions of Theorem~\ref{thm1} there exists $\eta>0$ such that for all $\lambda \in (r,s)$, all real $\delta$ with $|\delta| \le \frac{1}{8} \sqrt{(\lambda
-r)(s-\lambda)}$ and all subsets $A,B \subset \NN$ with $\max B < \min A$ it holds that
\begin{equation} \label{eq1.7}
\|\chi_A (J-(\lambda+i\delta))^{-1} \chi_B\| \le 4 \max\{(\lambda-r)^{-1}, (s-\lambda)^{-1}] \exp \left( -\eta \sqrt{(\lambda-r)(s-\lambda)} \sum_{k=\max B}^{\min A -1} \frac{1}{\lambda_k} \right).
\end{equation}
\end{theo}

Here $\chi_A$ and $\chi_B$ denote the multiplication operators with the characteristic functions of $A$ and $B$ and $\|\cdot\|$ the operator norm. Theorem~\ref{thm1} is a special case of Theorem~\ref{thm4}, where the energy is real ($\delta=0$), $A=\{n\}$ and $B=\{1\}$.

\vspace{.3cm}

\setcounter{equation}{0}
\section{Proofs of Theorems~\ref{thm1}, \ref{thm2} and \ref{thm3}} \label{sec3}

\subsection{Proofs of Theorem~\ref{thm1} and \ref{thm4}}


Let $\rho(n) := \sum_{k=1}^{n-1} 1/\lambda_k$ and, for $\gamma >0$
to be specified later, $\phi = e^{-\gamma \rho}$ as a
multiplication operator in $\ell^2(\NN)$. Then a calculation shows
that
\begin{equation} \label{eq2.1}
\phi^{-1} J \phi - J = A(\gamma) = \left( \begin{array}{cccc} 0 &
a_1 & & \\ b_1 & 0 & a_2 & \\ & b_2 & 0 & \ddots
\\ & & \ddots & \ddots \end{array} \right)
\end{equation}
is a non-symmetric Jacobi matrix with entries
\begin{equation} \label{eq2.2}
a_n = \lambda_n \left( \frac{\phi(n+1)}{\phi(n)} -1 \right), \quad
b_n = \lambda_n \left( \frac{\phi(n)}{\phi(n+1)} -1 \right), \quad
n\in \NN \;.
\end{equation}
To determine bounds on the real and imaginary part of $A$ we
verify by Taylor expansion that
\begin{eqnarray} \label{eq2.3}
a_n+b_n & = & \lambda_n \left( e^{\gamma (\rho(n+1)-\rho(n))} +
e^{-\gamma (\rho(n+1)-\rho(n))} - 2 \right) \nonumber \\ & = &
\frac{\gamma^2}{\lambda_n} + \Oo(\gamma^4/\lambda_n^3)
\end{eqnarray}
and
\begin{eqnarray} \label{eq2.4}
a_n - b_n & = & \lambda_n \left( e^{-\gamma (\rho(n+1)-\rho(n))} -
e^{\gamma (\rho(n+1) -\rho(n))} \right) \nonumber \\ & = &
-2\gamma + \Oo(\gamma^3/{\lambda_n}^2 ) \;.
\end{eqnarray}
We conclude that
\begin{equation} \label{eq2.5}
\mbox{\rm Re}\,A(\gamma) = \frac{1}{2} \left( \begin{array}{cccc}
0 & a_1+b_1 & & \\ a_1+b_1 & 0 & a_2+b_2 & \\ & a_2+b_2 & 0 &
\ddots
\\ & & \ddots & \ddots \end{array} \right) = \frac{\gamma^2}{2}
\left( \begin{array}{cccc}
0 & \lambda_1^{-1} & & \\ \lambda_1^{-1} & 0 & \lambda_2^{-1} & \\ & \lambda_2^{-1} & 0 & \ddots
\\ & & \ddots & \ddots \end{array} \right) + \Oo(\gamma^4)
\end{equation}
and
\begin{equation} \label{eq2.6}
\mbox{\rm Im}\,A(\gamma) = \frac{1}{2i} \left( \begin{array}{cccc}
0 & a_1-b_1 & & \\ b_1-a_1 & 0 & a_2-b_2 & \\ & b_2-a_2 & 0 &
\ddots
\\ & & \ddots & \ddots \end{array} \right) = i\gamma \left( \begin{array}{cccc}
0 & 1 & & \\ -1 & 0 & 1 & \\ & -1 & 0 & \ddots
\\ & & \ddots & \ddots \end{array} \right) + \Oo(\gamma^3)\;,
\end{equation}
where $\lambda_n \ge C>0$ uniformly in $n$ was used and error
terms refer to norm bounds. In fact, all we will use below are the
norm bounds $\| \mbox{\rm Re}\,A(\gamma)\| \le C_1 \gamma^2$ and
$\| \mbox{\rm Im}\,A(\gamma)\| \le C_2 \gamma$ with $C_1, C_2 $
depending only on $J$.

The following lemma provides the operator theoretic fact behind
the improvement of the Combes-Thomas method found in \cite{BCH}.
Extracting from the argument in \cite{BCH}, this lemma is stated
with proof in \cite[p.~60]{Stoll}.

\begin{lemma}
\label{lemma1}
Let $T$ be invertible self-adjoint with $d_+ = \mbox{\rm
dist}(0,\sigma(T) \cap (0,\infty))$ and \\ $d_- = \mbox{\rm
dist}(0,\sigma(T) \cap (-\infty,0))$. Let $S$ be self-adjoint,
$\|S\| \le 1$. Then for $\beta \in \RR$, $|\beta| \le \frac{1}{2}
\sqrt{d_+ \cdot d_-}$, the operator $T+i\beta S$ is invertible,
with
\begin{equation} \label{eq2.7}
\| (T+i\beta S)^{-1} \| \le 2 \cdot \max \left\{ \frac{1}{d_+},
\frac{1}{d_-} \right\} \;.
\end{equation}
\end{lemma}
\begin{rem}\label{remark1}
If $\sigma(T)\subset\RR_{+}$ (resp. $\sigma(T)\subset\RR_{-}$ )
then we put $ d_- = \infty $ ( resp. $ d_+ = \infty  $ ) .
\end{rem}

We apply Lemma~\ref{lemma1} to the operator $J-\lambda + A(\gamma) =
T + i\beta S$, with self-adjoint operators $T= J-\lambda + \mbox{\rm
Re}\,A(\gamma)$ and $S= \mbox{\rm Im}\,A(\gamma) / \| \mbox{\rm
Im}\,A(\gamma)\|$, and $\beta = \| \mbox{\rm Im}\,A(\gamma)\|$. As
$\| \mbox{\rm Re}\,A(\gamma)\| \le C_1\gamma^2$, the operator $T$
has a spectral gap $(r-\lambda +C_1\gamma^2, s-\lambda -C_1\gamma^2)
= (-d_-,d_+)$.

For a finite interval $(r,s)$ we choose $\gamma = \eta \sqrt{(\lambda
-r)(s-\lambda)}$, where
\begin{equation} \label{eq2.8}
\eta = \min \left\{ \frac{1}{4C_2}, \frac{1}{\sqrt{2C_1(s-r)}}
\right\} \;.
\end{equation}
Then $d_- \ge \frac{1}{2} (\lambda -r)$, $d_+ \ge \frac{1}{2}
(s-\lambda)$ and $|\beta| \le \frac{1}{4} \sqrt{(\lambda
-r)(s-\lambda)} \le \frac{1}{2} \sqrt{d_+ d_-}$. Applying
Lemma~\ref{lemma1} we know that $J-\lambda + A(\gamma)$ is
invertible and thus, using (\ref{eq2.7}),
\begin{eqnarray} \label{2.9}
\|\phi^{-1} (J-\lambda)^{-1} \phi \| & = & \| (J-\lambda +
A(\gamma))^{-1} \| \nonumber \\
& \le & 2 \max \left\{ \frac{1}{d_+}, \frac{1}{d_-} \right\}
\nonumber \\
& \le & 4 \max \left\{ \frac{1}{s-\lambda}, \frac{1}{\lambda-r}
\right\} \;.
\end{eqnarray}
We note that $|\langle (J-\lambda)^{-1} e_1, e_n \rangle | = \|
e_n (J-\lambda)^{-1} e_1\|$, where on the right we think of $e_n$
and $e_1$ as multiplication operators (and of $\|\cdot\|$ as the
operator norm). The proof of Theorem~\ref{thm1}(a) is thereby
completed through the estimate
\begin{eqnarray} \label{2.10}
\|e_n (J-\lambda)^{-1} e_1 \| & \le & \|e_n \phi\| \|\phi^{-1}
(J-\lambda)^{-1} \phi \| \|\phi^{-1}e_1\| \nonumber \\
& \le & e^{-\gamma \rho(n)} \cdot 4 \max \left\{
\frac{1}{s-\lambda}, \frac{1}{\lambda-r} \right\} \cdot 1 \;.
\end{eqnarray}

\vspace{.2cm}

Part (b) follows by simple modifications of the above argument. We are in the case $d_-=\infty$, meaning that there is no restriction on the size of $|\beta|$. One chooses $\gamma = \sqrt{(d-\lambda)/2C_1}$.

\vspace{.2cm}

The proof of Theorem~\ref{thm4} also follows the above lines with only minor changes. The imaginary part of the spectral parameter is included in $S$ by choosing $S= (\delta I+\mbox{\rm Im}\,A(\gamma)) / (|\delta|+\| \mbox{\rm
Im}\,A(\gamma)\|)$, and $\beta = |\delta|+\| \mbox{\rm Im}\,A(\gamma)\|$. The constant $\eta$ is modified to be the minimum of $\frac{1}{8C_2}$ and $\frac{1}{\sqrt{2C_1(s-r)}}$. The calculation (\ref{2.10}) applies in the same way to $\|\chi_A (J-(\lambda+i\delta))^{-1} \chi_B\|$ to give the bound (\ref{eq1.7}).

\vspace{.2cm}

Note that in the above argument crucial use was made of the fact
that the real part of $A(\gamma)$ is smaller than its imaginary
part, $\Oo(\gamma^2)$ as compared to $\Oo(\gamma)$. This is
exploited through Lemma~\ref{lemma1}, which can be interpreted as
saying that purely imaginary perturbations have a weaker effect on
the invertibility of self-adjoint operators than self-adjoint
perturbations. In the ``classical'' Combes-Thomas method only the
rougher bound $\|A(\gamma)\| = \Oo(\gamma)$ is used, yielding merely
linear dependence of the exponent in (\ref{eq1.4}) on the distance
of $\lambda$ to the spectrum.

\subsection{Proof of Theorem~\ref{thm2}}

We start with an improved version of Lemma~\ref{lemma1}.

\begin{lemma} \label{lemma2}
Let $T$ be invertible self-adjoint with positive $d_+ = \mbox{\rm
dist}(0,\sigma(T) \cap (0,\infty))$ and $d_- = \mbox{\rm
dist}(0,\sigma(T) \cap (-\infty,0))$. Let $S$ be self-adjoint,
$\|S\| \le 1$. Then for $\beta \in \RR$, $|\beta|<\sqrt{d_+ \cdot
d_-}$, the operator $T+i\beta S$ is invertible, with
\begin{equation} \label{eq2.11}
\| (T+i\beta S)^{-1} \| \le[\Delta_+ -(\Delta_+^2 +\beta^2 -
d_+\cdot d_-)^{1/2} ]^{-1}, \quad \Delta_+ =\frac{1}{2}(d_- + d_+).
\end{equation}
\end{lemma}
\quad
\begin{rem}\label{remark2}
If $\sigma(T)\subset\RR_{+}$ ( resp. $\sigma(T)\subset\RR_{-}$ )
then we take the limit of the right hand side of (\ref{eq2.11}) as
$d_{-}$ tends to $+\infty$ (resp. as $d_{+}$ tends to $ +\infty $
), and put no restriction on $\beta $ .
\end{rem}

\noindent {\em Proof:} Let us start with some geometrical ideas
which are behind the proof. Denote by $K$ be the upper half-disc
with the center at $(d_{+}-d_{-})/2$ ( on the x-axis ) of radius
$\Delta_{+}$. Then the length of the segment from $(0,0)$ to the
point of the intersection of the upper circle with the y-axis
equals $\sqrt{d_+ \cdot d_-}$ .  This explains why
$|\beta|<\sqrt{d_+ \cdot d_-}$ .

 Let $P_{+}$ (respectively $P_{-}$) be the
spectral projection of $T$ corresponding to the positive
(respectively negative) part of $\sigma(T).$ We have
$$ T = |T|(P_{+}-P_{-}) = T_{1}(d_{+}P_{+} - d_{-}P_{-}),$$
where $ T_{1} := |T|(\frac{P_{+}}{d_{+}} - \frac{P_{-}}{d_{-}}).$
Note that $ T_{1} $ commutes with $P_{+}$ and $P_{-}$ and
$T_{1}\geq{I}.$

Hence
\begin{eqnarray*}
\|(T+i\beta S)^{-1}\| & = & \|T_{1}^{-1/2}(d_{+}P_{+} - d_{-}P_{-}+
i\beta T_{1}^{-1/2} S T_{1}^{-1/2})^{-1}T_{1}^{-1/2}\| \\
& \leq & \|(d_{+}P_{+} - d_{-}P_{-}+ i\beta S_{1})^{-1}\|,
\end{eqnarray*}
with $S_{1}^* = S_{1}:= T_{1}^{-1/2} S T_{1}^{-1/2} $,
$\|S_{1}\|\leq{1}.$

Therefore the proof is reduced to the case $ T = d_{+}P_{+} -
d_{-}P_{-}.$ Denote by $\Delta_{-} = (d_{+} - d_{-})/2$ and $J_{1} =
P_{+} - P_{-}.$ Then
\begin{equation} \label{eq2.12}
 T+i\beta S = \Delta_{-}I +\Delta_{+}J_{1} + i\beta S = \Delta_{+}[J_{1}+\Delta_{+}^{-1}(\Delta_{-}I+i\beta S)]
\end{equation}
and
\begin{equation} \label{eq2.13}
\|\Delta_{-}I + i\beta S\|^2 = \|\Delta_{-}^2I+ \beta^2 S^*S\|
\leq{\Delta_{-}^2+ \beta^2} = \Delta_{+}^2 + (\beta^2 -
d_{-}d_{+}) .
\end{equation}

Using (\ref{eq2.12}) and (\ref{eq2.13}) we obtain
\begin{eqnarray*}
\|(T+i\beta S)^{-1}\| & \leq &
\Delta_{+}^{-1}\|[J_{1}+(\Delta_{-}I+i\beta
S)\Delta_{+}^{-1}]^{-1}\| \\
& \leq & \Delta_{+}^{-1}[1-(\Delta_{+}^2+\beta^2 -
d_-d_+)^{1/2}\Delta_{+}^{-1}]^{-1} \\ & = &
[\Delta_{+}-(\Delta_{+}^2+\beta^2 - d_-d_+)^{1/2}]^{-1}.
\end{eqnarray*}
This completes the proof of Lemma~\ref{lemma2}.

\vspace{.2cm}

Turning to the proof of Theorem~\ref{thm2} we fix
$\epsilon\in(0,1/2).$ Let
\begin{equation} \label{eq2.14}
\rho(n):=0 \:\mbox{for}\: n\le{N}, \: \rho(n) := \sum_{k=N}^{n-1}
1/\lambda_k \:\mbox{for} \: n>N,
\end{equation}
with $N$ to be chosen below. For $\gamma>0$ to be specified later,
let $\phi = e^{-\gamma \rho}$ be a multiplication operator in
$\ell^2(\NN)$. Then by repeating the calculations given in the proof
of Theorem~\ref{thm1} we find
\begin{equation} \label{eq2.15}
\phi^{-1} J \phi - J = A(\gamma) = \left( \begin{array}{cccc} 0 &
a_1 &  \\ b_1 & 0 & a_2 & \\ & b_2 & 0 & \ddots
\\ & & \ddots & \ddots \end{array} \right)
\end{equation}
is a non-symmetric Jacobi matrix with entries
\begin{equation} \label{eq2.16}
a_n = \lambda_n \left( \frac{\phi(n+1)}{\phi(n)} -1 \right), \quad
b_n = \lambda_n \left( \frac{\phi(n)}{\phi(n+1)} -1 \right), \quad
n\in \NN \;.
\end{equation}
We want to determine bounds on the real and imaginary part of
$A(\gamma)$. Note that
\begin{eqnarray}\label{eq2.17}
{\rm Re}\,A(\gamma) & = & SW + WS^{*}, \quad W = Diag (0,\dots,
\psi_1(N),\psi_1(N+1),\dots ), \\ \label{eq2.18} {\rm Im}\,A(\gamma) & =
& i(SZ-ZS^{*}), \quad Z = Diag(0,\dots,\psi_2(N),\psi_2(N+1),\dots),
\end{eqnarray}
where
\[ \psi_1(p)=
\sum_{k=1}^{\infty}\gamma^{2k}[(2k)!\lambda_p^{2k-1}]^{-1}, \quad
\psi_2(p)=-\sum_{k=1}^{\infty}\gamma^{2k-1}[(2k-1)!\lambda_p^{2k-2}]^{-1}.
\]

We apply Lemma~\ref{lemma2} to the operator $ J-\lambda + A(\gamma)
= T + i\beta S$, with self-adjoint operators $T = J-\lambda +
\mbox{\rm Re}\,A(\gamma)$ and $S= \mbox{\rm Im}\,A(\gamma) / \|
\mbox{\rm Im}\,A(\gamma)\|$, and $\beta = \| \mbox{\rm
Im}\,A(\gamma)\|$. \vspace{.1cm}
 Using (\ref{eq2.17}) and (\ref{eq2.18}) we have
\begin{eqnarray}\label{eq2.19}
\epsilon(N)=\| \mbox{\rm Re}\,A(\gamma)\|,
 \quad \delta(N)=\|
\mbox{\rm Im}\,A(\gamma) \|,
\end{eqnarray}
 where $\epsilon(N):= 2\sup_p{\psi_1(p)}$ and
$\delta(N):= 2\sup_p{(-\psi_2(p))}$.

Let $\gamma = (1/2 -
\epsilon)\sqrt{(\lambda-r)(s-\lambda)}$. Easy computations show
that
\begin{eqnarray}\label{eq2.20}
\epsilon(N)\leq{C_1(J)(\lambda-r)(s-\lambda)(\inf_{p\geq{N}}{\lambda_{p}})^{-1}}
\end{eqnarray} and
\begin{eqnarray}\label{eq2.21}
\delta(N)\leq{(1-2\epsilon)[(\lambda-r)(s-\lambda)]^{1/2}[1+\frac{1}{8}(s+r)^2(\inf{\lambda_{p}})^{-2}C_{2}(J)]}
\end{eqnarray}
for some positive constants $C_1(J)$ and $C_2(J)$ which depend
only on $J.$

 Note that the operator $T$ has a spectral gap $(r-\lambda
+\gamma^2\epsilon(N), s-\lambda -\gamma^2\epsilon(N)) =:
(-d_-,d_+)$. Denote $w(\lambda)=(\lambda-r)(s-\lambda).$ Evoking the
above formulas we obtain
\begin{eqnarray}\label{eq2.22}
d_{+}\geq{s-\lambda-C_1(J)w(\lambda)(\inf_{p\geq{N}}{\lambda_p})^{-1}}
\end{eqnarray}\,
and
\begin{eqnarray}\label{eq2.23}
d_{-}\geq{\lambda-r-C_1(J)w(\lambda)(\inf_{p\geq{N}}{\lambda_p})^{-1}}.
\end{eqnarray}
Let $r(N)= (\inf_{p\geq{N}}\lambda_p)^{-1}.$ Inequalities
(\ref{eq2.21}), (\ref{eq2.22}) and (\ref{eq2.23}) allow to compare
$\beta^2$ and $d_+\cdot d_- $. We find
\begin{eqnarray*}
d_+\cdot d_- -\beta^2 & \geq &
(s-\lambda-C_1(J)w(\lambda)r(N))(\lambda-r-C_1(J)w(\lambda)r(N)) \\
& & \mbox{}-(1-2\epsilon)^2w(\lambda)[1+\frac{1}{8}(s+r)^2
C_{1}(J)r(N)^2]^2 \\
& = & w(\lambda) \big[
(1-C_{1}(J)(\lambda-r)r(N))(1-C_1(J)(s-\lambda)r(N))
\\ & & \mbox{} -(1-2\epsilon)^2(1+\frac{1}{8}(s+r)^2 C_{1}(J)r(N)^2) \big] \\
& \geq & w(\lambda)[4\epsilon-4\epsilon^2 - C_{3}(J)r(N)],
\end{eqnarray*}
for some positive constant $C_{3}(J)$.

Choosing $N$ so large that $r(N)\leq4\epsilon^2(C_{3}(J))^{-1} $
the above inequalities show that
$$d_+\cdot d_-
-\beta^2\geq{4w(\lambda)\epsilon(1-2\epsilon)}\geq{\epsilon
w(\lambda)},$$ for $\epsilon\in (0,1/2)$.

Therefore Lemma \ref{lemma2} implies that $J-\lambda + A(\gamma)$ is
invertible and
\begin{eqnarray*}
\|\phi^{-1} (J-\lambda)^{-1} \phi \| & = & \| (J-\lambda +
A(\gamma))^{-1} \| \\
& \leq & [\Delta_+ -(\Delta_+^2 + \beta^2 - d_+\cdot d_- )^{1/2}
]^{-1} \\
& \leq & 2\Delta_+[(d_+\cdot d_- - \beta^2)]^{-1}
\\ & \leq & 2\Delta_+(\epsilon
w(\lambda))^{-1}=(s-r)(w(\lambda)\epsilon)^{-1}.
\end{eqnarray*}
Using the last inequalities and repeating the reasoning given at the
end of the proof of Theorem \ref{thm1} we get the desired estimate,
thus completing the proof of Theorem~\ref{thm2}

\begin{rem} \label{remark3}
{\rm The constant $\eta$ in Theorem~\ref{thm1} can be made
arbitrary small. In turn the choice of $\eta =1/2 - \epsilon $ in
Theorem~\ref{thm2} is optimal as will be shown below in
Example~\ref{example1}.}
\end{rem}

\subsection{Proof of Theorem~\ref{thm3}}

Fix $\lambda$ such that $\Re{\lambda}<d$ and $\epsilon \in (0,1)$.
Choose $N$ such that
\begin{eqnarray}\label{eq2.24}
\inf_{n\geq{N}}\lambda_n \geq1, \quad
2(\inf_{n\geq{N}}\lambda_n)^{-1}\exp[(1-\epsilon)^2(d
-\Re\lambda)]\leq\epsilon.
\end{eqnarray}
Fix an arbitrary $N_1>N$. Now define the sequence $\rho(n)$ by
\[ \rho(n)= \left\{
\begin{array}{ll} 0, & n\le N, \\ \sum_{k=N}^{n-1}(\sqrt{\lambda_k})^{-1}, &
N<n\leq{N_1}, \\ \sum_{k=N}^{N_1-1}(\sqrt{\lambda_k})^{-1}, & n>N_1.
\end{array} \right. \]
 Put $\gamma = (1-\epsilon) \sqrt{d-\Re{\lambda}}$ and consider the
 multiplication operator $\phi_{N_1} = e^{-\gamma \rho}$ in
$\ell^2(\NN)$. Note that $\phi_{N_1}$ is bounded and invertible (for
any $N_1>N $). By repeating the computation given in the proof of
Theorem~\ref{thm1} we have
\begin{equation}
\phi_{N_1}^{-1} J \phi_{N_1} - J =: A(\gamma)
\end{equation}
is a finite rank Jacobi operator with entries $a_k$ and $b_k$ (see
the proof of Theorem~\ref{thm1}), where
$a_1=\dots=a_{N-1}=b_1=\dots=b_{N-1}=0$ and $a_k=b_k =0$ for
$k\geq{N_{1}}$. By definition of $A(\gamma)$ one can easily check
that
\begin{eqnarray}\label{eq2.25}
\| \mbox{\rm Re}\,A(\gamma)\|\leq{\gamma^2
(1+2e^{\gamma^2}(\inf_{p\geq{N}}\lambda_p)^{-1})}
\end{eqnarray}

Define the operator $A $ in $\ell^2(\NN)$ on the domain $D(J)$
($D(A)=D(J)$) by
$$ Af:= Jf + \mbox{\rm Re}\,A(\gamma)f+i\mbox{\rm Im}\,A(\gamma)f .$$
 Hence for complex $\mu$ such that $\Re{\mu}<d-\|
\mbox{\rm Re}\,A(\gamma)\|$,  $A-\mu$ is invertible and
\begin{eqnarray}\label{eq2.26}
\|(A-\mu)^{-1}\| \leq[(d-\| \mbox{\rm
Re}\,A(\gamma)\|)-\Re\mu]^{-1}
\end{eqnarray}

Since
$$ A = J + A(\gamma) = \phi_{N_1}^{-1} J \phi_{N_1}$$
for $\lambda$ fixed as above we have
\begin{eqnarray}\label{eq2.27}
(A-\mu)^{-1} =\phi_{N_1}^{-1} (J-\mu)^{-1} \phi_{N_1}.
\end{eqnarray}
In turn (\ref{eq2.27}) implies that
$$|\langle (A-\lambda)^{-1}
e_1,e_n \rangle | = |\langle(J-\lambda)^{-1}\phi_{N_1}
e_1,\phi_{N_1}^{-1}e_n \rangle|=
\phi_{N_1}^{-1}(n)|\langle(J-\lambda)^{-1}e_1,e_n \rangle|,$$ for
any $N<n\leq{N_1}$.

 Using inequality (\ref{eq2.26}) we obtain
\begin{eqnarray}\label{eq2.28}
|\langle(J-\lambda)^{-1}e_1,e_n \rangle|\leq{\exp(-\gamma
\sum_{k=N}^{n-1} \frac{1}{\sqrt\lambda_k}) \:\frac{1}{ d-\|
\mbox{\rm Re}\,A(\gamma)\|)-\Re\lambda} }
\end{eqnarray}
provided that $N<n\leq{N_1}$ and $ \Re{\lambda}< d-\| \mbox{\rm
Re}\,A(\gamma)\|$.

Finally, using the definition of $\gamma$ as well as
(\ref{eq2.24}) and (\ref{eq2.25}) one can readily check that
\begin{eqnarray} \label{eq2.29}
\| \mbox{\rm Re}\,A(\gamma)\| \leq{(1-\epsilon)(d-\Re{\lambda})}.
\end{eqnarray}

Combining the above relations (\ref{eq2.28}) and (\ref{eq2.29}) and
the fact that $N_1$ was arbitrary completes the proof of
Theorem~\ref{thm3}.

\begin{rem} \label{remark4}

{\rm  Note that the above $N=N(\Re{\lambda},\epsilon)$ tends to
$\infty $ as either $\Re{\lambda}$ goes to $-\infty$ or
$\epsilon\longrightarrow 0 $. In the interesting region $\{\lambda :
\:\Re{\lambda}\approx{d} \} $ the dependence of $N$ on
$\Re{\lambda}$ disappears.}

\end{rem}

\setcounter{equation}{0}
\section{Optimality of the decay bounds}
\label{sec4}

In this section we discuss two specific models illustrating that the
bounds obtained above are optimal.

\begin{ex} \label{example1}

{\rm  The first model of unbounded Jacobi matrices $J$ concerns
Theorem \ref{thm2} and is given by
\begin{equation}\label{eq3.1}
q_n=0, \quad \lambda_n = n + c_n, \quad n\in \NN,
\end{equation}
where $(c_n)$ is a two-periodic sequence $(c_1, c_2, c_1, c_2,
\ldots)$ such that $c_1 \not= c_2$ and $\lambda_n
>0$ for all $n$.
Unbounded Jacobi matrices quite frequently fall into one of two
extreme cases, namely that either $\sigma_{ess}(J) =\RR$ or
$\sigma_{ess}(J) = \emptyset$. What makes the class (\ref{eq3.1})
interesting is that it is non-trivial in this respect. It can be
shown that $\sigma_{ess}(J) = \RR \setminus (-|\rho|, |\rho|)$,
where $\rho =c_1-c_2$. In fact, the spectrum of $J$ is purely
absolutely continuous in $\RR \setminus [-|\rho|, |\rho|]$. This
is proven in \cite{JM}and \cite{JNSJCAM} by finding asymptotics of
solutions of the equation $Ju = \lambda u $ and using the method
of subordinacy (for the a.c.\ spectrum) for this and related
classes of entries. In the same papers it is shown that
$\sigma(J)$ in $(-|\rho|, |\rho|)$ is empty. In other words we
have exactly the situation considered in Theorem \ref{thm2}, with
$ r=-|\rho| $ and $s=|\rho|$.

 Moreover, in \cite{JM} and \cite{JNSJCAM} it was shown that
for $\lambda\in(r,s)$ there exists a solution $u_1(n)$ of the
equation
\begin{equation}\label{eq3.2}
\lambda_{n-1} u(n-1) + q_n u(n) + \lambda_n u(n+1) = \lambda u(n),
\:n>1
\end{equation}
such that
\begin{equation}
[(J-\lambda)^{-1}e_1](n) = a u_{1}(n), \:n>1, \:a\in\mathbb{R}
\end{equation}
\begin{equation}\label{eq3.3}
\left(\begin{array}{c}u_{1}(2n-1)\\u_{1}(2n)\end{array}\right)=d_n\exp[-\sqrt{r^2-\lambda^2}\sum_{k=1}^{n}(2k)^{-1}]S(e_{-}+
o(1)),
\end{equation}
where $d_n= (-1)^n\prod_{k=1}^{n}(1-1/2k)$, $e_{-}=
\left(\begin{array}{c}0\\1\end{array}\right)$, and
$$ S : = \left( \begin{array}{cc}
1 & 1 \\
\frac{w+r}{\lambda}, & \frac{-w+r}{\lambda}
 \end{array} \right).$$

 Comparing (\ref{eq3.3}) with the estimate of Theorem~\ref{thm2}
 completes  the proof of the above mentioned sharpness of the
 result.}
 \end{ex}

\begin{rem} \label{remark5}

{\rm The above solution $u_1$ is the unique $\ell^2$-generalized
eigenfunction to which the bound found in Theorem~\ref{thm2}
applies. As $d_n \sim n^{-1/2}$ is subexponential, the decay of
$u_1$ at $+\infty$ is governed by $\exp \left(
-\sqrt{\rho^2-\lambda^2} \sum_{k=1}^n (2k)^{-1} \right)$. With $\eta
= 2^{-1}$, the $n$-dependence of the exponent through the factor
$\sum_{k=1}^n k^{-1} = \sum_{k=1}^n 1/\lambda_k$ and (for $\lambda$
close to the boundary of the essential spectrum) the dependence on
the distance of $\lambda$ to the spectrum given through
$\sqrt{\rho^2-\lambda^2} = \sqrt{(\rho - \lambda)(\lambda+\rho)}$
are exactly as found in Theorem~\ref{thm2}.}

\end{rem}

The next example we present below illustrates that the
 estimates of Theorem~\ref{thm3} are also sharp.

\begin{ex} \label{example2}
{\rm Consider the sequences given by
\begin{equation}\label{eq3.4}
q_n=-2n, \quad \lambda_n= n, \quad n\in \NN.
\end{equation}
Note that $J$ with entries defined by (\ref{eq3.4}) is bounded from
above by $-I$.

Fix $\lambda>-1 $ and $\epsilon\in(0,1).$ By the Birkhoff-Adams
theorem \cite{E} there exists a basis $u_{\pm}(n)$ of solutions of
(\ref{eq3.1}) with asymptotics given by
\begin{equation}\label{eq3.5}
u_{\pm}(n) = n^{-1/4}\exp(\pm2\sqrt{(\lambda+1)n})(1+o(1)).
\end{equation}

Applying Theorem~\ref{thm3} (in this case for $J$ bounded from above
by $(-I) $) there exists $N$ such that for $n>N$ we have

$$ |u_{-}(n)|\leq{(\epsilon(1+\lambda))^{-1}\exp[-(1-\epsilon)\sqrt{1+\lambda}\sum_{k=N}^{n}(\sqrt{k})^{-1}]}. $$

Since $\sum_{k=N}^{n}(\sqrt{k})^{-1}\simeq{2\sqrt{n}} $, by
comparing (\ref{eq3.5})  and the last estimate we conclude that the
value $1-\epsilon$ in the formula for $\gamma$ (see the the proof of
Theorem~\ref{thm3}) is optimal.}
\end{ex}

\noindent {\bf Problem.} In case that $J$ is a bounded Jacobi matrix
it is well known (see Theorem~2.3 in \cite{Shubin}) that the
spectrum of $J$ as an operator in $l^p$ does not depend on $p.$
However, for unbounded $J$ this result does not apply. Nevertheless,
in the case $ \sum_{n=1}^{\infty} \frac{1}{\lambda_n} = \infty $ the
estimates given in  Theorem~\ref{thm2} imply that $\sigma^2(J)
\supseteq{\sigma^p(J)}$ for any $p\in[1,\infty],$ where
$\sigma^p(J)$ is the spectrum of $J$ considered on the maximal
 domain in $l^p .$ This can be easily seen by applying the Schur test
 to $\langle (J-\lambda)^{-1} e_j, e_n \rangle$. Does the opposite inclusion also hold true?

\setcounter{equation}{0}
\section{A criterion for the existence of square-integrable solutions} \label{sec5}

While the results of this section could be stated and proven for general Jacobi matrices, we will for simplicity assume that $J$ is given by (\ref{eq1.1}) with zero-diagonal, $q_n=0$ for all $n\in\NN$.

We will compare $J$ with a second Jacobi matrix $J_0$, also with zero-diagonal and weights $\lambda_n^0$, $n\in \NN$. Both weight sequences satisfy (\ref{eq1.3}) to guarantee self-adjointness of $J$ and $J_0$.

We assume that the weights of $J$ and $J_0$ coincide on a sequence of increasing intervals: Let $x_k \in \NN$, $x_1 < x_2 < x_3 < \ldots$, $\ell_k \in \NN$ with $\lim_{k\to\infty} \ell_k =\infty$ and $x_k+ \ell_k < x_{k+1}-\ell_{k+1}$ for all $k$, and assume that
\begin{equation} \label{eq5.1}
\lambda_n = \lambda_n^0 \quad \mbox{for all $n\in \cup_k \{x_k-\ell_k-2, x_k+\ell_k+1\}$}.
\end{equation}
Denote $I_k := \{x_k-\ell_k, x_k+\ell_k\}$ and
\[
\Lambda_k := \max \{\lambda_n: x_k-\ell_k-2 \le n \le x_k+\ell_k+1\}.
\]

\begin{theo} \label{thm5.1}
Suppose that for all $\gamma>0$,
\begin{equation} \label{eq5.2}
\sum_{k=1}^{\infty} \Lambda_k (\Lambda_{k-1}+\Lambda_k+\Lambda_{k+1}) e^{-\gamma \ell_k/\Lambda_k} (x_{k+1}-x_{k-1}) < \infty,
\end{equation}
then for almost every $E\in \RR \setminus \sigma(J_0)$ there exists a non-trivial square-summable generalized eigenfunction of $J$ to $E$.
\end{theo}

Note that (\ref{eq5.2}) allows for situations where the distances between centers of barriers $x_{k+1}-x_k$ can grow significantly faster than their size $2 \ell_k$. This leads to applications where $J$ may contain a lot of spectrum outside the spectrum of $J_0$, including entire intervals. However, this spectrum can not be absolutely continuous and will typically be pure point, as follows from the following well-known general result.

\begin{proposi} \label{prop5.2}
Let $J$ be an unbounded self-adjoint Jacobi matrix in $\ell^2(\NN)$ and $I$ an open subset of $\RR$. Assume that for almost every $E\in I$ there exists a non-trivial square-summable generalized eigenfunction of $J$ to $E$. Then

(a) $\sigma_{ac}(J) \cap I = \emptyset$,

(b) $\sigma_c(J_{\lambda}) \cap I = \emptyset$ for almost every $\lambda \in \RR$, where $J_{\lambda} = J+ \lambda \langle \cdot, e_1 \rangle e_1$.
\end{proposi}

A proof of (b) for the case of discrete Schr\"odinger operators (i.e.\ all $\lambda_n=1$), using spectral averaging over the coupling parameter $\lambda$, is given in \cite{KMP1}. Due to the limit-point condition (\ref{eq1.3}) this proof extends to our setting. Part (a) follows from (b) as the absolutely continuous spectrum is invariant under rank one perturbations.

\vspace{.3cm}

The remainder of this section will be devoted to the proof of Theorem~\ref{thm5.1}. Note first that it suffices to show the existence of a non-trivial square-integrable generalized eigenfunction for almost every $E\in [\alpha,\beta]$, where the compact interval $[\alpha,\beta]$ is disjoint from $\sigma(J_0)$, as $\RR\setminus \sigma(J_0)$ can be exhausted by countably many such intervals.

For $E\in \RR$, $\eta>0$ and $z:=E+i\eta$ let
\begin{equation} \label{eq5.3}
u_{\eta}(n,E) := \langle (J-z)^{-1}e_1,e_n \rangle, \quad n\in \NN
\end{equation}
be the Weyl-solution, i.e.\ the unique $\ell_2$-solution of
\begin{equation} \label{eq5.4}
\lambda_{n-1}u_{n-1} + \lambda_n u_{n+1} = zu_n, \quad n\in \NN
\end{equation} satisfying the boundary condition $u_{\eta}(0,E)=-1$. The $m$-function is given by $m(E+i\eta) := u_{\eta}(1,E)$ and it is known that $m(E+i\eta) \to m(E)$ as $\eta\to 0$ exists and is non-zero for almost every $E\in \RR$. For arbitrary $z\in \CC$ let $\psi(\cdot,z)$ and $\varphi(\cdot,z)$ be the solutions of (\ref{eq5.4}) with $\psi(0,z)=-1$, $\psi(1,z)=0$, $\varphi(0,z)=0$ and $\varphi(1,z)=1$. Then
\[
u_{\eta}(n,E) = \psi(n,E+i\eta) + m(E+i\eta) \varphi(n,E+i\eta), \quad n\ge 0.
\]
We have pointwise in $n$ and $E$ that $\psi(n,E+i\eta) \to \psi(n,E)$ and $\varphi(n,E+i\eta) \to \varphi(n,E)$ as $\eta\to 0$. Thus, for almost every $E\in \RR$,
\begin{equation} \label{eq5.4a}
u(n,E):= \psi(n,E) +m(E)\varphi(n,E)
\end{equation}
exists and is a solution of (\ref{eq5.4}) with $z$ replaced by $E$. Theorem~\ref{thm5.1} is proven if we can show that, for a suitable $\eta_0>0$,
\begin{equation} \label{eq5.5}
\sup_{0<\eta\le \eta_0} \sum_{n=1}^{\infty} |u_{\eta}(n,E)|^2 <\infty \quad \mbox{for a.e.\ $E\in [\alpha,\beta]$},
\end{equation}
as this implies that $u(\cdot,E)$ is square-summable. One has $u(1,E)=m(E)$, thus $u(\cdot,E)$ is also non-trivial for almost every $E$.

In the following we will find bounds for the $\ell^2$-norm of $u_{\eta}$ separately within the barriers $I_k$ and for the intervals between any two given barriers.

To handle the barriers, define $\tilde{\ell}_k$ by rounding down $\ell_k/2$ to the closest integer and let $\tilde{I}_k := [x_k-\tilde{\ell}_k, x_k + \tilde{\ell}_k]$. Also, write $\chi_k$ and $\tilde{\chi}_k$ for the characteristic functions of $I_k$ and $\tilde{I}_k$. A calculation using (\ref{eq5.1}) shows that
\begin{equation} \label{eq5.6}
((J_0-(E+i\eta))\chi_k u_{\eta})(n) = \left\{ \begin{array}{ll} -\lambda_{n-1} u_{\eta}(n-1), & n\in \{x_k-\ell_k, x_k+\ell_k-1\}, \\ -\lambda_n u_{\eta}(n+1), & n\in \{x_k+\ell_k, x_k +\ell_k+1\}, \\ 0, & \mbox{else}. \end{array} \right.
\end{equation}
We write
\begin{equation} \label{eq5.7}
\tilde{\chi}_k u_{\eta} = \tilde{\chi_k} \chi_k u_{\eta} = \tilde{\chi}_k (J_0-(E+i\eta))^{-1} \chi_{U_k} (J_0-(E+i\eta)) \chi_k u_{\eta},
\end{equation}
where
\[
U_k := \{x_k-\ell_k-2, x_k-\ell_k-1, x_k-\ell_k, x_k+\ell_k, x_k+\ell_k+1, x_k+\ell_k+2\}.
\]
Plugging (\ref{eq5.6}) into (\ref{eq5.7}) yields
\begin{eqnarray} \label{eq5.8}
\sum_{n\in \tilde{I}_k} |u_{\eta}(n)|^2 & \le & 4\Lambda_k^2 \|\tilde{\chi}_k (J_0-(E+i\eta))^{-1} \chi_{U_k}\|^2 \sum_{n\in U_k} |u_{\eta}(n)|^2 \nonumber \\ & \le & 4\Lambda_k^2 a_k^2 \sum_{n\in U_k} |u_{\eta}(n)|^2,
\end{eqnarray}
where the elementary bound $(|a|+|b|+|c|+|d|)^2 \le 4(|a|^2+|b|^2+|c|^2+|d|^2)$ was used and we have set
\[
a_k := \sup_{0<\eta\le \eta_0, E\in [a,b]} \|\tilde{\chi}_k (J_0-(E+i\eta))^{-1} \chi_{U_k}\|.
\]

We now consider the intervals between barriers. For this let $J_k$ be the finite Jacobi matrix found by restricting $J$ to the interval $[x_k, x_{k+1}]$. We have, also considered after restricting to $[x_k,x_{k+1}]$,
\begin{equation} \label{eq5.9}
(1-\tilde{\chi}_k - \tilde{\chi}_{k+1}) u_{\eta} = (J_k-(E+i\eta))^{-1} (J_k-(E+i\eta))(1-\tilde{\chi}_k - \tilde{\chi}_{k+1}) u_{\eta}.
\end{equation}
Similar to (\ref{eq5.6}) we find
\begin{equation} \label{5.10}
((J_k-(E+i\eta))(1-\tilde{\chi}_k-\tilde{\chi}_{k+1})u_{\eta})(n) = \left\{ \begin{array}{ll} -\lambda_{n-1} u_{\eta}(n-1), & n\in \{x_k+\tilde{\ell}_k+1, x_k+\tilde{\ell}_k\}, \\ -\lambda_n u_{\eta}(n+1), & n\in \{ x_{k+1}-\tilde{\ell}_{k+1}-1, x_{k+1} -\tilde{\ell}_{k+1}\}, \\ 0, & \mbox{elsewhere in $[x_k,x_{k+1}]$}. \end{array} \right.
\end{equation}
With
\[ V_k := \{x_k+\tilde{\ell}_k-1, x_k+\tilde{\ell}_k, x_{k+1}-\tilde{\ell}_{k+1}, x_{k+1}-\tilde{\ell}_{k+1}+1\}
\]
we get from (\ref{eq5.9}) that
\begin{equation} \label{eq5.11}
\sum_{n=x_k+\tilde{\ell}_k+1}^{x_{k+1}-\tilde{\ell}_{k+1}-1} |u_{\eta}|^2 \le 4 (\Lambda_k^2 + \Lambda_{k+1}^2) \|(J_k-(E+i\eta))^{-1}\|^2 \sum_{n\in V_k} |u_{\eta}|^2.
\end{equation}
Given that $V_k \subset \tilde{I}_k \cup \tilde{I}_{k+1}$ we can bound the term $\sum_{n\in V_k} |u_{\eta}|^2$ on the right hand side of (\ref{eq5.11}) by two terms of the form (\ref{eq5.8}) to arrive at
\begin{equation} \label{eq5.12}
\sum_{n=x_k+\tilde{\ell}_k+1}^{x_{k+1}-\tilde{\ell}_{k+1}-1} |u_{\eta}|^2 \le \frac{16(\Lambda_k^2+\Lambda_{k+1}^2)}{\Delta_k(E)^2} \left( \Lambda_k^2 a_k^2 \sum_{n\in U_k} |u_{\eta}|^2 + \Lambda_{k+1}^2 a_{k+1}^2 \sum_{n\in U_{k+1}} |u_{\eta}|^2 \right).
\end{equation}
Here we have also used that
\[
\|(J_k-(E+i\eta))^{-1}\| \le \frac{1}{\Delta_k(E)},
\]
abbreviating $\Delta_k(E) := \mbox{dist}(E,\sigma(J_k))$.

Ultimately, we can bound the $\ell^2$-norm of $u_{\eta}$ over the entire interval $[x_k,x_{k+1}]$ by the sums in (\ref{eq5.8}) (for $\tilde{I}_k$ and $\tilde{I}_{k+1}$) and in (\ref{eq5.12}) to get
\begin{eqnarray} \label{eq5.13}
\sum_{n=x_k}^{x_{k+1}} |u_{\eta}(n)|^2 & \le & 4 \Lambda_k^2 a_k^2 \left( 1 + \frac{4(\Lambda_k^2+\Lambda_{k+1}^2)}{\Delta_k(E)^2} \right) \sum_{n\in U_k} |u_{\eta}(n)|^2 \nonumber \\
& & \mbox{}+ 4\Lambda_{k+1}^2 a_{k+1}^2 \left( 1 + \frac{4(\Lambda_k^2+\Lambda_{k+1}^2)}{\Delta_k(E)^2} \right) \sum_{n\in U_{k+1}} |u_{\eta}(n)|^2.
\end{eqnarray}
Summing over $k\ge K$ for any given $K\in \NN$ yields
\begin{equation} \label{eq5.14}
\sum_{n=x_K}^{\infty} |u_{\eta}(n)|^2 \le 16 \sum_{k=K}^{\infty} b_k(E) \sum_{n=x_k-\ell_k-2}^{x_k+\ell_k+2} |u_{\eta}(n)|^2,
\end{equation}
where
\begin{equation} \label{eqbk}
b_k(E) := \Lambda_k^2a_k^2\left( 1+\frac{\Lambda_{k-1}^2+\Lambda_k^2}{\Delta_k(E)^2} + \frac{\Lambda_k^2+\Lambda_{k-1}^2}{\Delta_k(E)^2} \right).
\end{equation}
We will now use
\begin{lemma} \label{lem5.2}
Under the conditions of Theorem~\ref{thm5.1} it holds that
\begin{equation} \label{eq:bto0}
\lim_{k\to\infty} b_k(E) =0 \quad \mbox{for a.e.\ $E\in [\alpha,\beta]$}.
\end{equation}
\end{lemma}
Before proving Lemma~\ref{lem5.2}, we show how it is used to complete the proof of Theorem~\ref{thm5.1}. Let $E\in [\alpha,\beta]$ be such that $b_k(E)\to 0$ and also $m(E+i\eta)\to m(E)\in \CC$ as $\eta\to 0$. By Lemma~\ref{lem5.2} and our remarks at the beginning of the proof, this holds for almost every $E\in [\alpha,\beta]$.

There exists $k_0(E)$ such that $b_k(E) \le 1/32$ for all $k\ge k_0(E)$. By our general assumptions, any given $n\in \NN$ is not covered by more than two of the intervals $[x_k-\ell_k-2,x_k+\ell_k+2]$. Thus, if we pick $K=k_0(E)$ in (\ref{eq5.14}), it follows that
\begin{eqnarray*}
\sum_{n\in\NN} |u_{\eta}(n)|^2 & = & \sum_{n=1}^{x_K-1} |u_{\eta}(n)|^2 + \sum_{n=x_K}^{\infty} |u_{\eta}(n)|^2 \nonumber \\
& \le & \sum_{n=1}^{x_K-1} |u_{\eta}(n)|^2 +\frac{1}{2} \sum_{n\in\ZZ} |u_{\eta}(n)|^2,
\end{eqnarray*}
or
\begin{equation} \label{eqfinish}
\sum_{n\in\NN} |u_{\eta}(n)|^2 \le 2 \sum_{n=1}^{x_K-1} |u_{\eta}(n)|^2.
\end{equation}
The crucial fact is that this holds uniformly in $\eta$. By (\ref{eq5.4a}), $\lim u_{\eta}(n) =u(n)$ exists for each of the finitely many $n \in [1,x_K-1]$. Thus the right hand side of (\ref{eqfinish}) is uniformly bounded in $\eta\in (0,1]$. This concludes the proof of (\ref{eq5.5}) and thus of Theorem~\ref{thm5.1}.

\vspace{.3cm}

We finally \emph{prove Lemma~\ref{lem5.2}}, where Theorem~\ref{thm4} is the crucial ingredient. The interval $[\alpha,\beta]$ has positive distance $d_0 := \,\mbox{dist}([\alpha,\beta], \sigma(J_0))$ from the spectrum of $J_0$. At this point we choose $\eta_0 := d_0/8$, meaning that $\eta_0 \le \sqrt{(\lambda-r)(s-\lambda)}/8$ for all $\lambda \in [\alpha,\beta]$, where $r:= \max(\sigma(J_0) \cap (-\infty, \alpha))$ and $s:= \min(\sigma(J_0) \cap (\beta,\infty))$. Thus, by Theorem~\ref{thm4}, there exist constants $C<\infty$ and $\gamma_0>0$ such that
\[
a_k \le C \exp \left( -\gamma_0 \sum_{n=x_k-\ell_k-1}^{x_k-\tilde{\ell}_k-1} \frac{1}{\lambda_n} \right) + C \exp \left( -\gamma_0 \sum_{n=x_k+\tilde{\ell}_k}^{x_k+\ell_k} \frac{1}{\lambda_n} \right).
\]
Here we have split the four points in $U_k$ into the two pairs to the right and left of $\tilde{I}_k$ and applied (\ref{eq1.7}) separately. Using the definition of $\Lambda_k$ and that $\tilde{\ell}_k\le \ell_k/2$, we arrive at
\begin{equation} \label{eq5.15}
a_k \le 2C e^{-\eta_0 \ell_k/2\Lambda_k}.
\end{equation}
Let $A_k:= \{E\in [\alpha,\beta]: \Delta_k(E)<\alpha_k\}$ with $\alpha_k$ to be chosen later. $J_k$ has $x_{k+1}-x_k+1$ eigenvalues. Thus $|A_k| \le 2\alpha_k (x_{k+1}-x_k+1) \le 4\alpha_k(x_{k+1}-x_k)$. Suppose that
\begin{equation} \label{eq5.16}
\sum_k \alpha_k (x_{k+1}-x_k) < \infty.
\end{equation}
Then $\sum_k |A_k| < \infty$ and it follows from the Borel-Cantelli lemma that for almost every $E\in [\alpha,\beta]$ there exists $k_0(E) \in \NN$ such that $\Delta_k(E)\ge \alpha_k$ for $k\ge k_0(E)$. For such $E$ and $k$ we conclude from (\ref{eqbk}) and (\ref{eq5.15}) that
\begin{equation} \label{eq5.17}
b_k(E) \le C\Lambda_k^2 e^{-\gamma0 \ell_k/2\Lambda_k} \left(1+\frac{\Lambda_{k-1}^2+\Lambda_k^2}{\alpha_{k-1}^2} + \frac{\Lambda_k^2+\Lambda_{k+1}^2}{\alpha_k^2} \right).
\end{equation}
Let $\gamma_1 \in (0,\gamma_0/2)$ and
\[
\alpha_k^2 := \Lambda_k^2 (\Lambda_k^2+\Lambda_{k+1}^2) e^{-\gamma_1 \ell_k/\Lambda_k} + \Lambda_{k+1}^2 (\Lambda_k^2+\Lambda_{k+1}^2) e^{-\gamma_1 \ell_{k+1}/\Lambda_{k+1}}.
\]
With this choice of $\alpha_k$, (\ref{eq5.2}) shows the summability of $|A_k|$. Moreover, from (\ref{eq5.17}) we get
\[ b_k(E) \le C\left(\Lambda_k^2 e^{-\gamma_0 \ell_k/2\Lambda_k} + e^{-(\frac{1}{2}\gamma_0-\gamma_1) \ell_k/\Lambda_k} \right),
\]
which tends to $0$ as $k\to\infty$ by (\ref{eq5.2}). Lemma~\ref{lem5.2} is proven.

\setcounter{equation}{0}
\section{A class of Jacobi matrices with a mobility edge} \label{sec6}

As an application of Theorem~\ref{thm5.1} we now provide a class of unbounded Jacobi matrices, which exhibit a transition from spectral regions with purely absolutely continuous spectrum to a region with dense pure point spectrum. Due to the connection of spectral and transport properties this is called a mobility edge in the physics literature. Our example is motivated by classes of bounded Jacobi matrices considered in \cite{Stolz2} which exhibit the same behavior. Here, as well as in \cite{Stolz2}, a purely absolutely continuous Jacobi matrix is subjected to a slowly oscillating perturbation, which generates regions of dense pure point spectrum. As opposed to the examples in \cite{Stolz2}, in our example an additional periodic modulation sequence is needed to open up a gap in the purely absolutely continuous spectrum of the unperturbed unbounded Jacobi matrix with weights $n^{\alpha}$, $0<\alpha<1$.

Define the weights by
\begin{equation} \label{eq6.1}
\lambda_n = n^{\alpha} + c_n \varphi(n^{\gamma}).
\end{equation}
Here $\varphi:\RR \to \RR$ is twice continuously differentiable and periodic, i.e.\ $\varphi(x+T)=\varphi(x)$ for some $T>0$ and all $x\in \RR$. Furthermore $0\le \varphi(x)\le 1$ for all $x$, $\inf \varphi =0$ and $\sup \varphi =1$. The sequence $c_n$, $n\in\NN$, is $2$-periodic, $c_n=c_1$ for all odd $n$ and $c_n=c_2$ for all even $n$ with $c_1>0$, $c_2>0$ and $c_1\not= c_2$. Write $c:= |c_1-c_2|$.

\begin{theo} \label{thm6.1}
Let $J$ be the Jacobi matrix on $\ell^2(\NN)$ with zero-diagonal and weights given by (\ref{eq6.1}). Suppose that $0<\alpha<1$ and $0<\gamma<(1-\alpha)/2$. Then

(a) $\sigma(J) = \RR$,

(b) $J$ is purely absolutely continuous in $\RR \setminus [-c,c]$,

(c) $J$ has no absolutely continuous spectrum in $(-c,c)$,

(d) $J_{\lambda} := J+\lambda \langle \cdot, e_1\rangle e_1$ has pure point spectrum in $(-c,c)$ for almost every $\lambda \in \RR$.
\end{theo}

We start with the proof of part (a). Denote by $J_0$ the Jacobi operator with zero diagonal and the weights $\lambda_n^{(0)} = n^{\alpha}$, $\alpha \in (0,1]$. Let $(J_0 u)_n = \lambda u_n$, $n>1$, $\lambda \in \RR$. Using Theorem~3.2 of \cite{JNJFA} we know that
\begin{equation} \label{eq51a1}
|u_n| \sim n^{-\alpha/2} \quad \mbox{as $n\to\infty$}.
\end{equation}
Let $\varepsilon_i := i^{\frac{\alpha+\gamma-1}{\gamma}+\delta}$, where $0<\delta<\frac{1-\alpha-\gamma}{\gamma}$. Assume that $\varphi(x_0)=0$ and consider the sequence $\{x_i\}\subset \RR$ such that $x_i^{\gamma} =x_0+iT$. It follows that $x_i \sim T^{1/\gamma} i^{1/\gamma}$.

For $n_i := [x_i]$ choose a sequence $\{\Delta_i\} \subset \NN$ of even numbers such that
\begin{equation} \label{eq51a2}
\Delta_i \ge \varepsilon_i n_i^{1-\gamma}
\end{equation}
and
\begin{equation} \label{eq51a3}
\Delta_i \le M n_i^{1-\gamma-\varepsilon}
\end{equation}
for some $M>0$ and $0<\varepsilon < 1-\alpha-\gamma -\gamma\delta$. Consider the sequence of intervals $I_i := [n_i, n_i+\Delta_i]$. Denote $\tilde{n}_i := n_i + \Delta_i/2$. We let $\beta_i :=2/\Delta_i$ and claim that the sequence $\{v^{(i)}\} \in \ell^2$ given by
\[ v^{(i)}(n) := \left\{ \begin{array}{ll} 0, & n\not\in I_i \\ u_n [1+\beta_i(n-\tilde{n}_i)], & n_i \le n \le \tilde{n}_i \\ u_n[1-\beta_i(n-\tilde{n}_i)], & \tilde{n}_i < n \le n_i+\Delta_i, \end{array} \right. \]
is a Weyl sequence for $J$ at the point $\lambda$:

By definition of $v^{(i)}$ one can verify that
\begin{equation} \label{eq51a4}
\|v^{(i)}\|_2 \sim (\Delta_i n_i^{-\alpha})^{1/2}.
\end{equation}

Moreover, observe that the above choice of $\Delta_i$ implies
\[
\sup_{n\in I_i} |\varphi(n^{\gamma})| \to 0 \quad \mbox{as $i\to\infty$}.
\]
In fact
\begin{eqnarray*}
|\varphi(n^{\gamma})| & = & |\varphi(n^{\gamma})-\varphi(x_i^{\gamma})| \le \|\varphi'\|_{\infty} |(n_i+\Delta_i)^{\gamma}-n_i^{\gamma}| \\
& \le & C n_i^{\gamma} \frac{\Delta_i}{n_i} \le CM n_i^{-\varepsilon}
\end{eqnarray*}
for some $C>0$ and $i$ sufficiently large. Here we have used (\ref{eq51a3}). Hence
\begin{equation} \label{eq51a5}
\frac{\|(J-\lambda)v^{(i)}\|}{\|v^{(i)}\|} \le \frac{\|(J_0-\lambda)v^{(i)}\|}{\|v^{(i)}\|} + M_1 n_i^{-\varepsilon}
\end{equation}
for some $M_1>0$ and large $i$.

We find for $n_i< n < \tilde{n_i}$ (the "$+$"-signs in the following calculation) and for $\tilde{n}_i < n < n_i +\Delta_i$ (the "$-$"-signs), respectively,
\begin{eqnarray*}
[(J_0-\lambda)v^{(i)}](n) & = & [\pm (n-1)^{\alpha} u_{n-1}(n-1-\tilde{n}_i) \pm n^{\alpha} u_{n+1}(n+1-\tilde{n}_i) -\lambda (\pm u_n) (n-\tilde{n}_i)] \beta_i \\
& = & \pm \beta_i (n-\tilde{n}_i) [ (n-1)^{\alpha} u_{n-1} +n^{\alpha} u_{n+1}-\lambda u_n] \\
& & \pm [(n-1)^{\alpha} \beta_i (-u_{n-1}) +n^{\alpha} u_{n+1} \beta_i -\lambda \beta_i u_n].
\end{eqnarray*}
Note that the first term in the last expression vanishes as $J_0 u -\lambda u=0$. The only other values of $n$ such that $[(J_0-\lambda)v^{(i)}](n)\not= 0$, $n= n_i$, $\tilde{n}_i$, $n_i+\Delta_i$, give slightly different expressions which do not contribute significantly to $\|(J_0-\lambda)v^{(i)}\|$. All this and the bound (\ref{eq51a1}) lead to
\begin{equation} \label{eq51a6}
\|(J_0-\lambda)v^{(i)}\|^2 \le M_2 (n_i+\Delta_i)^{2\alpha} \beta_i^2 n_i^{-\alpha} (\Delta_i+3).
\end{equation}
Combining (\ref{eq51a4}), (\ref{eq51a5}) and (\ref{eq51a6}) we find
\[
\frac{\|(J-\lambda)v^{(i)}\|}{\|v^{(i)}\|} \le M_3 \beta_i n_i^{\alpha} + M_1 n_1^{-\varepsilon}
\]
for some $M_3>0$ and large $i$. From (\ref{eq51a2}) it follows that $\beta_i n_i^{\alpha} \le M_4 i^{-\delta} \to 0$ as $i\to\infty$. This completes the proof of part (a).

\vspace{.3cm}

(b) Fix $\lambda \in \{x\in \RR: |x|>c\}$. We want to study the behavior of generalized eigenfunctions of $Ju=\lambda u$. Denote by $B_n$ the transfer matrix of $J$ given by
\[
B_n := \left( \begin{array}{cc} 0 & 1 \\ -\frac{\lambda_{n-1}}{\lambda_n} & \frac{\lambda}{\lambda_n} \end{array} \right).
\]
As in our previous works \cite{JNKrein, JNJFA, JNSJCAM} it is of advantage to compute the products $B_{2n} B_{2n-1}$. We have
\begin{eqnarray} \label{eq51a7}
\frac{\lambda_{2n-1}}{\lambda_{2n}} & = & 1+ (2n)^{-\alpha} \frac{c_1 \varphi((2n-1)^{\gamma}) -c_2 \varphi((2n)^{\gamma})-\alpha (2n)^{\alpha-1}}{1+c_2 \varphi((2n)^{\gamma})(2n)^{-\alpha}} +p_n, \\
\label{eq51a8}
\frac{\lambda_{2n-2}}{\lambda_{2n-1}} & = & 1+ (2n)^{-\alpha} \frac{c_2 \varphi((2n-2)^{\gamma}) -c_1 \varphi((2n-1)^{\gamma})-\alpha (2n)^{\alpha-1}}{(1-1/2n)^{\alpha}+c_1 \varphi((2n-1)^{\gamma})(2n)^{-\alpha}}+z_n,
\end{eqnarray}
where $p_n$, $z_n \in l^1$. Since $0<\gamma <(1-\alpha)/2$ and $\varphi'$ and $\varphi''$ are bounded, the sequences $\{\varphi((2n-2)^{\gamma})\}$, $\{\varphi((2n-1)^{\gamma})\}$ and $\{\varphi((2n)^{\gamma})\}$ belong to the class $D^2((2n)^{\alpha})$, see Example~2.1 in \cite{J2M2}. Recall here that for a non-negative weight sequence $\mu=\{\mu_n\}$ the class $D^2(\mu)$ is defined in \cite{J2M2} as $\{x\in l^{\infty}: \Delta x \in l^2(\mu), \Delta^2 x \in l^1(\mu)\}$, where a sequence $x=\{x_n\}$ is in $l^p(\mu)$ if $\sum_n |x(n)|^p \mu(n) <\infty$ and $\Delta$ is the forward difference operator, i.e.\ $(\Delta x)_n = x_{n+1}-x_n$. It is clear that $\{n^{-\alpha}\}$ is also in $D^2((2n)^{\alpha})$.

The class $D^2((2n)^{\alpha}$ is closed under multiplication and division (by a sequence separated from zero), see Lemma~2.2 in \cite{J2M2}. Therefore we can write (\ref{eq51a7}) and (\ref{eq51a8}) as
\begin{eqnarray} \label{eq51a9}
\frac{\lambda_{2n-1}}{\lambda_{2n}} = 1 + (2n)^{-\alpha} [(c_1-c_2)\varphi((2n-1)^{\gamma})+x_n] +r_n, \\
\label{eq51a10}
\frac{\lambda_{2n-2}}{\lambda_{2n-1}} = 1+ (2n)^{-\alpha} [(c_2-c_1) \varphi((2n-2)^{\gamma})+y_n] +s_n,
\end{eqnarray}
where $\{x_n\}, \{y_n\}$ belong to $D^2((2n)^{\alpha})$, $\lim x_n = \lim y_n =0$, and $\{r_n\}$, $\{s_n\} \in l^1$.

By the same arguments
\begin{equation} \label{eq51a11}
\frac{1}{\lambda_{2n-1}} = (2n)^{-\alpha} t_n, \quad \frac{1}{\lambda_{2n}} = (2n)^{-\alpha} w_n,
\end{equation}
where $\{t_n\}$, $\{w_n\} \in D^2((2n)^{\alpha})$, $\lim t_n = \lim w_n = 1$.

Using (\ref{eq51a9}), (\ref{eq51a10}) and (\ref{eq51a11}) we obtain
\[
B_{2n} B_{2n-1} = -I + (2n)^{-\alpha} V(n) + R(n),
\]
with $\{\|R(n)\|\} \in l^1$. The matrix $V(n)$ has entries
\begin{eqnarray*}
V_{11}(n) & = & -(c_2-c_1) \varphi((2n-2)^{\gamma})-y_n, \\
V_{12}(n) & = & \lambda t_n, \\
V_{21}(n) & = & -\lambda w_n \{1+(2n)^{-\alpha} [(c_2-c_1) \varphi(2n-2)^{\gamma} +y_n]\}, \\
V_{22}(n) & = & -[(c_1-c_2) \varphi((2n-1)^{\gamma}) +x_n] +\lambda^2 (2n)^{-\alpha} t_n w_n.
\end{eqnarray*}
Thus $V(n)$ is in $D^2((2n)^{\alpha})$.

Note that
\begin{eqnarray*}
\mbox{discr}\,V(n) & := & (\mbox{tr}\,V(n))^2 - 4 \det V(n) \\
& = & -4 [\lambda^2 + (c_2-c_1)^2 \varphi((2n-1)^{\gamma}) \varphi((2n)^{\gamma})] + o(1),
\end{eqnarray*}
as $n\to\infty$. It follows that $\limsup \,\mbox{discr}\,V(n) <0$. Thus we have verified all the assumptions to apply Theorem~5.1 of \cite{J2M2}, and in particular the asymptotic formula (5.17) there. This shows that $Ju=\lambda u$ has no subordinated solutions in the sense of \cite{Khan/Pearson}. It follows from the results of \cite{Khan/Pearson} that $J$ is purely absolutely continuous in $\RR \setminus [-c,c]$.

\vspace{.3cm}

To prove parts (c) and (d), let $J_0$ be the Jacobi matrix on $\ell^2(\NN)$ with weights $\lambda_n^0 = n^{\alpha}+c_n$. It was shown in \cite{JNSJCAM} that $J_0$ has at most finitely many spectral points in $(-c,c)$. Thus $(-c,c) \setminus \sigma(J_0)$ is the union of finitely many open intervals $U_j$.

Fix a $U_j$ and a compact interval $[\alpha,\beta] \subset U_j$. We will show that $J$ has a non-trivial square-summable generalized eigenfunction for almost every $\lambda \in [\alpha,\beta]$. This implies that the assumptions of Proposition~\ref{prop5.2} hold with $I=(-c,c)$ as $(-c,c) \setminus \sigma(J_0)$ can be exhausted by countably many such intervals. Thus (c) and (d) follow.

Choose $x_0 \in [0,T)$ such that $\varphi(x_0)=1$ and let $\delta := \,\mbox{dist}([\alpha,\beta],\pm c)$. There exists $\varepsilon>0$ such that $\varphi(x) \ge 1-\delta/(2\max \{c_1,c_2\})$ for all $x \in [x_0-\varepsilon, x_0+\varepsilon]$. For every $k\in \NN$, let $x_k$ be the integer closest to the center of the interval $[(x_0+kT-\varepsilon)^{1/\gamma}, (x_0+kT+\varepsilon)^{1/\gamma}]$ and $\ell_k$ an integer approximately equal to a quarter of the length of this interval. Let $J_{\varepsilon}$ be the Jacobi matrix with weights
\[ \lambda_n^{\varepsilon} := \left\{ \begin{array}{ll} \lambda_n, & n \in \cup_k [x_k-\ell_k-2,x_k+\ell_k+1], \\ n^{\alpha}+c_n, & \mbox{else}. \end{array} \right. \]
Then $|\lambda_n^{\varepsilon}-(n^{\alpha}+c_n)|\le \delta/2$ for all $n$ and thus $[\alpha,\beta]$ is contained in a spectral gap of $J_{\varepsilon}$. We will complete our proof by showing that the assumptions of Theorem~\ref{thm5.1} are fulfilled with $J_{\varepsilon}$ in place of $J_0$.

For this, note that $x_k \sim (k_0+kT)^{1/\gamma}$ and that $x_{k+1}-x_k \sim \frac{T}{\gamma} (x_0+kT)^{1/\gamma -1}$ and $\Lambda_k \sim (x_0+kT)^{\alpha/\gamma}$ are polynomially bounded in $k$. On the other hand, $\ell_k \sim \frac{2\varepsilon}{\gamma} (x_0+kT)^{1/\gamma -1}$ and thus
\[
\ell_k/\Lambda_k \sim \frac{2\varepsilon}{\gamma} (x_0+kT)^{\frac{1-\alpha}{\gamma}-1}.
\]
The latter is polynomially growing in $k$ since $1-\alpha > (1-\alpha)/2>\gamma$. This implies that (\ref{eq5.2}) holds for all $\eta>0$.

\vspace{.3cm}

\noindent {\bf Acknowledgements:} We would like to thank B.\ Mityagin for making us aware of the references \cite{Demko} and \cite{Kershaw}. J.J.\ is supported in part by MSHE
grant N N201 426533, S.N.\ supported in part by ''INTAS'' and in
part by RFBR grant 06-01-00249, G.S.\ supported in part by NSF grant
DMS-0653374


{\small Institute of Mathematics, Polish Academy of Sciences, ul.\
sw.\ Tomasza 30, 31-027 Krakow, Poland, najanas@cyf-kr.edu.pl}

{\small Department of Mathematical Physics, Institute of Physics,
St.\ Petersburg University, Ulianovkaia 1, 198904 St.\ Petergoff,
St.\ Petersburg, Russia, naboko@snoopy.phys.spbu.ru}

{\small Department of Mathematics, University of Alabama at
Birmingham, CH 452, Birmingham, AL 35294, USA, stolz@math.uab.edu}


\begin{thebibliography}{99}

\bibitem{BCH} J.~M.~Barbaroux, J.-M.~Combes and P.~D.~Hislop:
Localization near band edges for random Schr\"{o}dinger operators,
{\it Helv.~Phys.~Acta}~{\bf 70} (1997), 16-43.

\bibitem{Breuer} Jonathan Breuer: Spectral and dynamical
properties of certain random Jacobi matrices with growing
parameters, Preprint http://front.math.ucdavis.edu/0708.0670

\bibitem{Berezanskii} Yu.~M.~Berezanskii: Expansions in
eigenfunctions of selfadjoint operators, {\it Translations of
Mathematical Monographs}, Vol.~{\bf 17}, Amer.~Math.~Soc.,
Providence, 1968

\bibitem{CT} J.-M.~Combes and L.~Thomas: Asymptotic behaviour of
eigenfunctions for multiparticle Schr\"{o}dinger operators, {\it
Comm.~Math.~Phys.}~{\bf 34} (1973), 251-270.

\bibitem{Demko} S.~Demko: Inverses of band matrices and local convergence of spline projections, {\it Siam~J.~Numer.~Anal.}~{\bf 14} (1977), 616-619.

\bibitem{E}
S.~N.~Elaydi:  An Introduction to Difference Equations, Springer,
New York, Inc., 1999.

\bibitem{GJMS} Y.~A.~Gordon, V.~Jaksic, S.~Molchanov and B.~Simon: Spectral Properties of Random Schr\"odinger Operators with Unbounded Potentials, {\it Commun.~Math.~Phys.}~{\bf 157} (1993), 23-50.

\bibitem{JM} J.~Janas and M.~Moszy\'{n}ski: Asymptotics of solutions
of difference equations and spectral properties of some Jacobi
matrices, {\it J.~Approx.~Th.}~{\bf 120} (2003), 309-336.

\bibitem{J2M2} J.~Janas and M.~Moszy\'{n}ski: New discrete Levinson type asymptotics of solutions of linear systems, {\it Jour.~Difference~Eq.~Appl.}~{\bf 12} (2006), 133-163.

\bibitem{JNKrein} J.~Janas and S.~Naboko: Asymptotics of
generalized eigenvectors for unbounded Jacobi matrices with
power-like weights, Pauli matrices, commutation relations and
Cesaro averaging, {\it Oper.~Th.~Advan.~Appl.}, Vol.~{\bf 117}
(2000), 165-186.

\bibitem{JNJFA} J.~Janas and S.~Naboko: Jacobi matrices with power
like weigths --- grouping in blocks approach, {\it
J.~Funct.~Anal.}~{\bf 166} (1999), 218-243.

\bibitem{JNSJCAM} J.~Janas, S.~Naboko and G.Stolz: Spectral theory for a class
of periodically perturbed unbounded Jacobi matrices:elementary
methods, {\it J.~Comp.~Appl.Math.}~{\bf 171} (2004), 265-276.

\bibitem{Kershaw} D.~Kershaw: Inequalities on the elements of the inverse of a certain tridiagonal matrix, {\it Math.~Comput.}~{\bf 24} (1970), 155-158.

\bibitem{Khan/Pearson} S.~Khan and D.~B.~Pearson: Subordinacy and spectral theory for infinite matrices, {\it Helv.~Phys.~Acta}~{\bf 65} (1992), 505-527.

\bibitem{KMP1} W.~Kirsch, S.~Molchanov and L.~Pastur: The one-dimensional Schr\"odinger operator with unbounded potential: the pure point spectrum, {\it Funktsional.~Anal.~i~Prilozhen.}~{\bf 24} (1990), 14-25 (in Russian), translation in {\it Funct.~Anal.~Appl.}~{\bf 24} (1990), 176-186 (1991).

\bibitem{KMP2} W.~Kirsch, S.~Molchanov, L.~Pastur: One-dimensional Schr\"odinger operators with high potential barriers, pp.~163-170, in: Oper.~Theory~Adv.~Appl.~{\bf 57}, Birkh\"auser, Basel, 1992.

\bibitem{KMPV} W.\ Kirsch, S.\ Molchanov, L.\ Pastur and B.\ Vainberg: Quasi 1D localization: deterministic and random potentials, {\it Markov Process.~Related Fields}~{\bf 9} (2003), 687-708.

\bibitem{Shubin} M.A.~Shubin: Pseudodifference operators and their
Green functions, {\it Sib.~Math.~Zh.}, Vol. ~{\bf 49} (1985),
652-671 (in Russian).

\bibitem{Stoll} P.~Stollmann: Caught by disorder: bound states
in random media, {\it Progress in Mathematical Physics}, Vol.~{\bf
20}, Birkh\"{a}user, Boston, 2001.

\bibitem{Stolz} G.~Stolz: Localization for Schr\"odinger operators with effective barriers, {\it J.~Funct.~Anal.}~{\bf 146} (1997), 416-429.

\bibitem{Stolz2} G.~Stolz: Spectral theory for slowly oscillating potentials. I.~Jacobi matrices, {\it Manuscripta Math.}~{\bf 84} (1994), 245-260.
\end{thebibliography}
\end{document}